\documentclass{amsart}
\usepackage{amssymb,amsmath,amsthm}
\usepackage[mathscr]{euscript}
\usepackage[latin1]{inputenc}

\usepackage{graphicx}
\input{epsf}

\usepackage[bbgreekl]{mathbbol}
\DeclareSymbolFontAlphabet{\mathbb}{AMSb}
\DeclareSymbolFontAlphabet{\mathbbb}{bbold}

\newcounter{sec}

\newcounter{punct}[sec]

\def\punct{\refstepcounter{punct}{\arabic{sec}.\arabic{punct}.  }}

\newtheorem{theorem}{Theorem}[sec]
\newtheorem{proposition}[theorem]{Proposition}

\newtheorem{lemma}[theorem]{Lemma}

\newtheorem{corollary}[theorem]{Corollary}

\def\COUNTERS{\addtocounter{sec}{1}
              \setcounter{punct}{0}
          \setcounter{equation}{0}
          \setcounter{theorem}{0}
            \setcounter{figure}{0}
          }
          
          \def\sm{\smallskip}

\begin{document}

\newcommand{\supp}{\mathop {\mathrm {supp}}\nolimits}
\newcommand{\rk}{\mathop {\mathrm {rk}}\nolimits}
\newcommand{\Aut}{\mathop {\mathrm {Aut}}\nolimits}
\newcommand{\Out}{\mathop {\mathrm {Out}}\nolimits}
\renewcommand{\Re}{\mathop {\mathrm {Re}}\nolimits}
\newcommand{\Inn}{\mathop {\mathrm {Inn}}\nolimits}
\newcommand{\Char}{\mathop {\mathrm {Char}}\nolimits}
\newcommand{\ch}{\cosh}
\newcommand{\sh}{\sinh}
\newcommand{\Sp}{\mathop {\mathrm {Sp}}\nolimits}
\newcommand{\SOS}{\mathop {\mathrm {SO^*}}\nolimits}
\newcommand{\Ams}{\mathop {\mathrm {Ams}}\nolimits}
\newcommand{\Gms}{\mathop {\mathrm {Gms}}\nolimits}

\def\0{\mathbf 0}

\def\ov{\overline}
\def\wh{\widehat}
\def\wt{\widetilde}
\def\pol{\twoheadrightarrow}

\newcommand{\dd}[1]{\,d\,{\overline{\overline{#1}}} }

\renewcommand{\rk}{\mathop {\mathrm {rk}}\nolimits}
\renewcommand{\Aut}{\mathop {\mathrm {Aut}}\nolimits}
\newcommand{\Ob}{\mathop {\mathrm {Ob}}\nolimits}
\renewcommand{\Re}{\mathop {\mathrm {Re}}\nolimits}
\renewcommand{\Im}{\mathop {\mathrm {Im}}\nolimits}
\newcommand{\sgn}{\mathop {\mathrm {sgn}}\nolimits}
\newcommand{\ver}{\mathop {\mathrm {vert}}\nolimits}
\newcommand{\val}{\mathop {\mathrm {val}}\nolimits}
\newcommand{\edge}{\mathop {\mathrm {edge}}\nolimits}
\newcommand{\germ}{\mathop {\mathrm {germ}}\nolimits}
\newcommand{\PB}{{\mathop {\EuScript {PB}}\nolimits}}

\def\bfa{\mathbf a}
\def\bfb{\mathbf b}
\def\bfc{\mathbf c}
\def\bfd{\mathbf d}
\def\bfe{\mathbf e}
\def\bff{\mathbf f}
\def\bfg{\mathbf g}
\def\bfh{\mathbf h}
\def\bfi{\mathbf i}
\def\bfj{\mathbf j}
\def\bfk{\mathbf k}
\def\bfl{\mathbf l}
\def\bfm{\mathbf m}
\def\bfn{\mathbf n}
\def\bfo{\mathbf o}
\def\bfp{\mathbf p}
\def\bfq{\mathbf q}
\def\bfr{\mathbf r}
\def\bfs{\mathbf s}
\def\bft{\mathbf t}
\def\bfu{\mathbf u}
\def\bfv{\mathbf v}
\def\bfw{\mathbf w}
\def\bfx{\mathbf x}
\def\bfy{\mathbf y}
\def\bfz{\mathbf z}

\def\bfA{\mathbf A}
\def\bfB{\mathbf B}
\def\bfC{\mathbf C}
\def\bfD{\mathbf D}
\def\bfE{\mathbf E}
\def\bfF{\mathbf F}
\def\bfG{\mathbf G}
\def\bfH{\mathbf H}
\def\bfI{\mathbf I}
\def\bfJ{\mathbf J}
\def\bfK{\mathbf K}
\def\bfL{\mathbf L}
\def\bfM{\mathbf M}
\def\bfN{\mathbf N}
\def\bfO{\mathbf O}
\def\bfP{\mathbf P}
\def\bfQ{\mathbf Q}
\def\bfR{\mathbf R}
\def\bfS{\mathbf S}
\def\bfT{\mathbf T}
\def\bfU{\mathbf U}
\def\bfV{\mathbf V}
\def\bfW{\mathbf W}
\def\bfX{\mathbf X}
\def\bfY{\mathbf Y}
\def\bfZ{\mathbf Z}

\def\frA{\mathfrak A}
\def\frB{\mathfrak B}
\def\frD{\mathfrak D}
\def\frS{\mathfrak S}
\def\frL{\mathfrak L}
\def\frP{\mathfrak P}
\def\frQ{\mathfrak Q}
\def\frR{\mathfrak R}
\def\frT{\mathfrak T}
\def\frG{\mathfrak G}
\def\frg{\mathfrak g}
\def\frh{\mathfrak h}
\def\frf{\mathfrak f}
\def\frk{\mathfrak k}
\def\frl{\mathfrak l}
\def\frm{\mathfrak m}
\def\frn{\mathfrak n}
\def\fro{\mathfrak o}
\def\frp{\mathfrak p}
\def\frq{\mathfrak q}
\def\frr{\mathfrak r}
\def\frs{\mathfrak s}
\def\frt{\mathfrak t}
\def\fru{\mathfrak u}
\def\frv{\mathfrak v}
\def\frw{\mathfrak w}
\def\frx{\mathfrak x}
\def\fry{\mathfrak y}
\def\frz{\mathfrak z}

\def\bfw{\mathbf w}
%%% END MATHBF
%%%%%%%%%%%%%%%%%%%%%%%%%%%%%%%
%%%%%%%%%%%%%%%%%%%%%%%%%%%%%%%%%
%%% BEGIN MATHBB

\def\R {{\mathbb {R} }}
 \def\C {{\mathbb C }}
  \def\Z{{\mathbb Z}}
  \def\H{{\mathbb H}}
\def\K{{\mathbb K}}
\def\N{{\mathbb N}}
\def\Q{{\mathbb Q}}
\def\A{{\mathbb A}}
\def\O {{\mathbb O }}
\def\G {{\mathbb G }}

\def\T{\mathbb T}
\def\P{\mathbb P}
\def\SS{\mathbb S}

\def\G{\mathbb I}

\def\cD{\EuScript D}
\def\cL{\EuScript L}
\def\cK{\EuScript K}
\def\cM{\EuScript M}
\def\cN{\EuScript N}
\def\cP{\EuScript P}
\def\cQ{\EuScript Q}
\def\cR{\EuScript R}
\def\cT{\EuScript T}
\def\cW{\EuScript W}
\def\cY{\EuScript Y}
\def\cF{\EuScript F}
\def\cG{\EuScript G}
\def\cZ{\EuScript Z}
\def\cI{\EuScript I}
\def\cB{\EuScript B}
\def\cA{\EuScript A}
\def\cE{\EuScript E}
\def\cC{\EuScript C}

\def\bbA{\mathbb A}
\def\bbB{\mathbb B}
\def\bbD{\mathbb D}
\def\bbE{\mathbb E}
\def\bbF{\mathbb F}
\def\bbG{\mathbb G}
\def\bbI{\mathbb I}
\def\bbJ{\mathbb J}
\def\bbL{\mathbb L}
\def\bbM{\mathbb M}
\def\bbN{\mathbb N}
\def\bbO{\mathbb O}
\def\bbP{\mathbb P}
\def\bbQ{\mathbb Q}
\def\bbS{\mathbb S}
\def\bbT{\mathbb T}
\def\bbU{\mathbb U}
\def\bbV{\mathbb V}
\def\bbW{\mathbb W}
\def\bbX{\mathbb X}
\def\bbY{\mathbb Y}

\def\U{\mathbb{U}}

\def\kappa{\varkappa}
\def\epsilon{\varepsilon}
\def\phi{\varphi}
\def\le{\leqslant}
\def\ge{\geqslant}

\def\B{\mathrm B}

\def\la{\langle}
\def\ra{\rangle}
\def\tri{\triangleright}

\def\lambdA{{\boldsymbol{\lambda}}}
\def\alphA{{\boldsymbol{\alpha}}}
\def\betA{{\boldsymbol{\beta}}}
\def\mU{{\boldsymbol{\mu}}}
\def\nU{{\boldsymbol{\nu}}}
\def\omegA{{\boldsymbol{\omega}}}
\def\GammA{{\mathbbb{\Gamma}}}
\def\taU{{\boldsymbol{\tau}}}
\def\thetA{{\boldsymbol{\theta}}}
\def\sigmA{{\boldsymbol{\sigma}}}
\def\kappA{{\boldsymbol{\varkappa}}}

\def\const{\mathrm{const}}
\def\rem{\mathrm{rem}}
\def\even{\mathrm{even}}
\def\SO{\mathrm{SO}}
\def\SL{\mathrm{SL}}
\def\PSL{\mathrm{PSL}}
\def\cont{\mathrm{cont}}
\def\Isom{\mathrm{Isom}}
\def\Isolated{\mathrm{Isolated}}
\def\coun{\mathrm{count}}
\def\S{\mathrm{S}}
\def\out{\mathrm{out}}

\def\Un{\operatorname{U}}
\def\GL{\operatorname{GL}}
\def\Mat{\operatorname{Mat}}
\def\End{\operatorname{End}}
\def\Mor{\operatorname{Mor}}
\def\Aut{\operatorname{Aut}}
\def\inv{\operatorname{inv}}
\def\red{\operatorname{red}}
\def\Ind{\operatorname{Ind}}
\def\dom{\operatorname{dom}}
\def\im{\operatorname{im}}
\def\md{\operatorname{mod\,}}
\def\indef{\operatorname{indef}}
\def\Gr{\operatorname{Gr}}
\def\Pol{\operatorname{Pol}}
\def\UU{\operatorname{U}}

\def\arr{\rightrightarrows}
\def\bs{\backslash}

\def\cH{\EuScript{H}}
\def\cO{\EuScript{O}}
\def\cQ{\EuScript{Q}}
\def\cL{\EuScript{L}}
\def\cX{\EuScript{X}}
\def\cJ{\EuScript{J}}

\def\Di{\Diamond}
\def\di{\diamond}

\def\fin{\mathrm{fin}}
\def\ThetA{\boldsymbol {\Theta}}

\def\0{\boldsymbol{0}}

\def\F{\,{\vphantom{F}}_2F_1}
\def\FF{\,{\vphantom{F}}_3F_2}
\def\H{\,\vphantom{H}^{\phantom{\star}}_2 H_2^\star}
\def\HH{\,\vphantom{H}^{\phantom{\star}}_3 H_3^\star}
\def\Ho{\,\vphantom{H}_2 H_2}

\def\1{{\boldsymbol{1}}}
\def\2{{\boldsymbol{2}}}

\newcommand{\res}{\mathop{\mathrm{res}}\nolimits}

\begin{center}

\bf \Large
On the Dotsenko-Fateev complex twin 
\\
of the Selberg integral and its extensions

\bigskip

{\sc Yury A. Neretin}
\footnote{Supported by the grant FWF--Austrian Scintific Funds, P31591.}

\end{center}

{\small The  Selberg integral has a twin (`the Dotsenko--Fateev integral')
of the following form.  We replace real variables $x_k$ in the integrand 
$\prod |x_k|^{\sigma-1}\,|1-x_k|^{\tau-1} \prod|x_k-x_l|^{2\theta}$ of the Selberg
integral by complex variables $z_k$, 
integration over a cube we  replace by an integration over the whole
complex space $\mathbb{C}^n$. According to Dotsenko, Fateev, and  Aomoto, such integral is a product of Gamma functions.
We define and evaluate a family of beta integrals over spaces $\mathbb{C}^m\times \mathbb{C}^{m+1}\times \dots \times \mathbb{C}^n$, which for $m=n$ gives the complex twin of
the Selberg integral mentioned above  (with
 three additional integer parameters).

}

\section{Results of the paper}

\COUNTERS

Recall that the Selberg integral (Selberg, 1944) is given by
\begin{multline}
S_n(\sigma,\tau;\theta):=\int_0^1\dots \int_0^1\prod_{j=1}^n |x_k|^{\sigma-1}|1-x_k|^{\tau-1}
\prod_{1\le k<l\le n}|x_k-x_l|^{2\theta}\prod_{j=1}^n dx_j
=\\=
n!\prod_{k=1}^n\frac{\Gamma(\sigma+(k-1)\theta)\,\Gamma(\tau+(k-1)\theta)\,\Gamma(k\theta)}{\Gamma(\sigma+\tau+(n+k-2)\theta)\,\Gamma(\theta)},
\label{eq:selberg}
\end{multline}
see a detailed discussion in Andrews, Askey, Roy \cite{AAR}, Chapter 8, for original Selberg's proof, see Lique, Thibon \cite{LT}.
Dotsenko and Fateev (1985, \cite{DF}, formula (B.9))  published (without a proof)  the following identity, which looks like a  twin of the Selberg integral, 
\begin{multline}
S_n^\C(\sigma,\tau;\theta):=
\\:=
\int_{\C^n} \prod_{k=1}^n  |z_k|^{2\sigma-2}|1-z_k|^{2\tau-2} \prod_{1\le k<l\le n}|z_k-z_l|^{4\theta} \prod_k d\Re z_k\,d\Im z_k
=\\=
\frac 1{n!}\prod_{j=1}^n\frac{ \sin\pi(\sigma+(j-1)\gamma)\,\sin\pi(\tau+(j-1)\theta)\,\sin \pi(j\theta)}{\sin\pi(\sigma+\tau+(n+j-2)\theta)\,\sin\pi(\theta)} 
\cdot
S(\sigma,\tau;\theta)^2.
\label{eq:df}
\end{multline}
This integral was rediscovered by 
Aomoto \cite{Aom} (who also presented a proof), see also Mimachi, Yoshida \cite{MY} and  a review of integrals of the Selberg type integrals by Forrester, Warnaar  \cite{For}.

The domain of convergence of this integral is narrow, see \cite{For},
\begin{align}
&\Re\sigma>0,\quad\Re\tau>0, \quad \Re \theta >-\tfrac  1n, 
\label{eq:ds-convergence-1}
\\
&\Re(\sigma+ (n-1)\theta)>0, \quad \Re(\tau+ (n-1)\theta)>0;
\label{eq:ds-convergence-2}
\\
&\Re(\sigma+\tau+(n-1)\theta)<1,\quad \Re(\sigma+\tau+2(n-1)\theta)<1.
\label{eq:ds-convergence-3}
\end{align}

\sm

{\bf \punct Some notation.}
For a complex variable $z$ we denote
$$
\dd{z}= d\Re z\,\,d\Im z=\tfrac 1{2i}dz\, d\ov z.
$$
By
$$
\bfa=a\bigl|a'
$$
we denote a pair of complex numbers such that $a-a'$ is integer. 
Denote by $\Lambda$ the set of such pairs.
It  splits into a countable  family of complex lines $a-a'=k$, where $k$ ranges in $\Z$.
Denote  
$$\1:=1\bigl|1,\qquad \2:=2\bigl|2, \qquad \lfloor\bfa\rfloor=\tfrac12 \Re(a+a').$$
For a complex $z$ 
denote
$$
z^\bfa=z^a\,\ov z^{a'}=|z|^{a+a'} e^{i(\arg z)(a-a')}.
$$
In particular, 
$$z^\1=|z|^2,\qquad (-1)^\bfa=(-1)^{a-a'}.$$

\sm

{\bf \punct Gamma-function of the complex field.}
Following Gelfand, Graev, and Retakh \cite{GGR}, we define the Gamma function of the complex field by
\begin{multline}
\Gamma^\C(\bfa)=\Gamma^\C(a|a'):= \frac 1\pi
\int_\C z^{\bfa-\1} e^{2i\Re z} \dd{z}:=\\=
i^{a-a'}
\frac{\Gamma(a)}{\Gamma(1-a')}=
i^{a'-a} \frac{\Gamma(a')}{\Gamma(1-a)}=
\frac{i^{a'-a}}\pi \Gamma(a)\,\Gamma(a')\sin \pi a'
.
\label{eq:gamma}
\end{multline}
The integral conditionally converges if $0<\lfloor\bfa\rfloor<1$, it admits a meromorphic continuation
to each  complex line $a-a'=k$. 

The {\it beta function of the complex field} is defined by
\begin{equation}
\B^\C(\bfa,\bfb):=\frac 1\pi \int_\C z^{\bfa-\1} (1-z)^{\bfb-\1}\,\dd{z}=
\frac{\Gamma^\C(\bfa)\,\Gamma^\C(\bfb)} {\Gamma^\C(\bfa+\bfb)}.
\label{eq:beta}
\end{equation}
The domain of (absolute) convergence is 
$$
\lfloor\bfa\rfloor>0,\quad\lfloor\bfb\rfloor>0,\quad \lfloor \bfa+\bfb\rfloor<1.
$$

\sm

{\bf \punct Results of the paper.} Denote by $\cZ$ a collection of point configurations
$$
z_{11}\in \C,\, (z_{21},z_{22})\in \C^2, \dots, (z_{n1},\dots,\,z_{nn})\in\C^n.
$$
Denote by $\cC_n$ the set of all possible $\cZ$, 
so  $\cC_n\simeq\C^{n(n+1)/2}$.
Denote
$$
\dd{\cZ}:=\prod_{j=1}^n \prod_{\alpha=1}^j \dd{z}_{j\alpha}.
$$
Let $\sigmA_1$, \dots, $\sigmA_n$, $\taU_1$, \dots, $\taU_n$, $\thetA_{11}$, $\{\thetA_{2j}\}_{j=1,2}$,
 \dots, $\{\thetA_{(n-1)j}\}_{j=1,2,\dots,n-1}\in \Lambda$.
 
 \begin{theorem}
 \label{th:1}
The following identity holds
\begin{multline}
\int\limits_{\vphantom{\bigl|}\cC_{\,n}}
 \prod_{j=1}^{n-1}\prod_{\alpha=1}^j z_{j\alpha}^{\sigmA_j-\sigmA_{j+1}-\thetA_{j\alpha}}(1-z_{j\alpha})^{\taU_j-\taU_{j+1}-\thetA_{j\alpha}}
 \prod_{p=1}^n z_{np}^{\sigmA_n-\1}(1-z_{np})^{\taU_n-\1}\times
 \\
\times \prod_{j=1}^{n-1} \frac{\prod\limits_{1\le \alpha\le j;\, 1\le p\le j+1}(z_{j\alpha}-z_{(j+1)p})^{\thetA_{j\alpha}-\1}}
 {\prod\limits_{1\le\alpha,\beta\le j;\,\alpha\ne\beta}(z_{j\alpha}-z_{j\beta})^{\thetA_{j\alpha}-\1}}  \prod_{1\le p<q\le n}(z_{np}-z_{nq})^{\1}\,\,
 \dd{\cZ}
 =
 \\=
 (-1)^{\sum_{j=1}^{n-1}\sum_{\alpha=1}^j \thetA_{j\alpha}}\pi^{n(n+1)/2}
 \prod_{j=1}^{n} j!
 \times \\ \times
  \prod_{1\le\alpha\le j\le n-1} \Gamma^\C(\thetA_{j\alpha})
 \prod_{j=1}^n
 \frac{\Gamma^\C(\sigmA_j)\,\Gamma^\C(\taU_j)}{\Gamma^\C(\sigmA_j+\taU_j+\sum_{1\le\alpha\le j-1}\thetA_{(j-1)\alpha})}.
 \label{eq:main}
\end{multline}
The domain of convergence of the integral is
$$
\lfloor\sigmA_j\rfloor>0,\quad \lfloor\taU_j\rfloor>0, 
\quad \lfloor\thetA_{j\alpha}\rfloor>0,
\quad \lfloor\sigmA_j+\taU_j+\sum_{1\le\alpha\le j-1}\thetA_{(j-1)\alpha}\rfloor< 1.
$$
\end{theorem}

This integral admits a further extension.
Denote by $\cC_n^m$
the space $\C^m\times\C^{m+1}\times \dots \times \C^n$ for $m\le n$.
We denote points of this space
by 
$$
\cZ:=\bigl((z_{m1}, \dots, z_{mm}), \dots, (z_{n1}, \dots, z_{nn})).
$$
We denote
$
\dd{\cZ}:=\prod_{j=m}^n \prod_{\alpha=1}^j \dd{z}_{j\alpha}
$.
Consider the following collection of parameters $\in\Lambda$:
\begin{align*}
&\wh\sigmA_n^m=\{\sigmA_m, \dots,\sigmA_n\}=:\{\sigmA_j\}_{j=m,\dots, n}, 
\qquad
\wh \taU^m_n=\{\taU_j\}_{j=m,\dots, n},
\\
& \wh\thetA_n^m=\{\thetA_{j\alpha}\}_{j=m,\dots,n-1;\,\alpha=1,\dots,j},\qquad \nU.
\end{align*}

\begin{theorem}
\label{th:2}
Denote
\begin{multline}
J_n^m(\wh\sigmA_n^m,\wh\taU_n^m;\nU,\wh\thetA_n^m)=
\\
=
\int_{\cC_{\,n}^{\,\vphantom{q_B}m}} \prod_{j=m}^{n-1}\prod_{\alpha=1}^j z_{j\alpha}^{\sigmA_j-\sigmA_{j+1}-\thetA_{j\alpha}}(1-z_{j\alpha})^{\taU_j-\taU_{j+1}-\thetA_{j\alpha}}
\cdot \prod_{p=1}^n z_{np}^{\sigmA_n-\1}(1-z_{np})^{\taU_n-\1}
\times\\\times
\prod_{1\le\alpha,\beta\le m;\, \alpha\ne \beta} (z_{m\beta}-z_{m\alpha})^{\nU-\frac\1\2}
\times\\\times \prod\limits_{j=m}^{n-1}
\frac{ \prod\limits_{1\le\alpha\le j;\, 1\le p\le j+1} (z_{j\alpha}-z_{(j+1)p})^{\thetA_{j\alpha}-\1} }
{  \prod\limits_{1\le \alpha,\beta\le j;\,\alpha\ne\beta}(z_{j\alpha}-z_{j\beta})^{\thetA_{j\alpha}-\1} }
\cdot \prod_{1\le p<q\le n} (z_{np}-z_{nq})^\1\,\dd{\cZ}.
\label{eq:trapezoid-0}
\end{multline}
The following identity holds:
\begin{multline}
J_n^m(\wh\sigmA_n^m,\wh\taU_n^m;\nU,\wh\thetA_n^m)=\\=
(-1)^{\nU m(m-1)/2+ \sum_{j=m}^{n-1}\sum_{\alpha=1}^j \thetA_{j\alpha}}\,\,
\pi^{n(n+1)/2-m(m-1)/2} \prod_{j=m}^n j!\,
\times \\ \times % \Gamma^\C(\nU)^{-m+1}  
 \prod_{j=1}^{m-1} \frac{\Gamma^\C(\sigmA_m+(m-j)\nU)\, \Gamma^\C(\taU_m+(m-j)\nU)\,\Gamma^\C((j+1)\nU)}
 {\Gamma^\C(\sigmA_m+\taU_m+(2m-j-1)\nU) \,\Gamma^\C(\nU)
 }
  \times\\ 
 \times
  \prod_{j=m}^{n-1}
 \prod_{\alpha=1}^j
%  \prod_{m\le\alpha\le j\le n-1}
  \Gamma^\C(\thetA_{j\alpha})
 \prod_{j=m}^n
 \frac{\Gamma^\C(\sigmA_{j})\,\Gamma^\C(\taU_{j})}{\Gamma^\C(\sigmA_j+\taU_j+\sum_{\alpha=1}^{j-1}\thetA_{(j-1)\alpha})}.
 \label{eq:trapezoid-1}
\end{multline}
The following conditions are sufficient for absolute convergence of the integral
\begin{align}
&\lfloor\sigmA_j\rfloor>0,\, \lfloor\taU_j\rfloor>0, 
 \, \lfloor\thetA_{j\alpha}\rfloor>0,
 \, \lfloor\sigmA_j+\taU_j+\sum_{1\le\alpha\le j-1}\thetA_{(j-1)\alpha}\rfloor< 1\,\,\text{for $j> m$};
 \label{eq:nec}
\\
&\lfloor\nU \rfloor>0, \,\,  \lfloor m\nU \rfloor<1,\,\, 
\lfloor \sigmA_m+\taU_m+2(m-1)\nU\rfloor<1,
\end{align}
the conditions \eqref{eq:nec} also are necessary.
\end{theorem}

{\sc Remark.} The factor in the second line of the right hand side
of \eqref{eq:trapezoid-0} has a similar factor in the denominator of the third line (for $j=m$), we can unite them,
but the formula is more readable without this combination. \hfill $\boxtimes$

\sm 

Next, we set $m=n$ in \eqref{eq:trapezoid-0}. Then two factor $\prod_{j=m}^{n-1}$  of the integrand disappear. We denote $\sigmA:=\sigmA_n$, $\taU:=\taU_n$,
$\thetA:=\nU$, 
$z_j=z_{nj}$ (other parameters in our case are absent) 
and  come to the following extended Dotsenko--Fateev--Aomoto integral:

\begin{corollary}
\begin{multline}
\int\limits_{\C^n}\prod_{j=1}^n z_j^{\sigmA-\1} (1-z_j)^{\taU-\1}
\prod_{1\le i<j\le n}(z_i-z_j)^{2\thetA}\prod_{j=1}^n \dd{z}_j
=\\=(-1)^{\thetA \cdot n(n-1)/2}
n!\,\pi^n\,\prod_{j=1}^n\frac{\Gamma^\C(\sigmA+(j-1)\thetA)\,\Gamma^\C(\taU+(j-1)\thetA)\,\Gamma^\C(j\thetA)}{\Gamma^\C(\sigmA+\taU+(n+j-2)\thetA)\,\Gamma^\C(\thetA)}.
\label{eq:dfa}
\end{multline}
\end{corollary}

Formally, this formula modulo a constant factor  can be obtained from the Selberg integral \eqref{eq:selberg} by
a substitution
\begin{equation}
x\mapsto z,\quad
\sigma\mapsto\sigmA,\quad \tau\mapsto\taU,\quad \theta\mapsto\thetA,\quad 1\to\1,\quad \Gamma\mapsto \Gamma^\C.
\label{eq:substitution}
\end{equation}
We also get 3 additional integer parameters in comparison to
\eqref{eq:df}, namely, $\sigma-\sigma'$, $\tau-\tau'$, $\theta-\theta'$.

\sm

{\sc Remarks.} a)
Our way of derivation is valid if we add the condition $\lfloor \thetA\rfloor>0$ to 
\eqref{eq:ds-convergence-1}--\eqref{eq:ds-convergence-3}, it is easy to omit it if we know the final formula.

\sm

b) The paper by Dotsenko and Fateev \cite{DF} contains another integral (B.10) over $\C^m\times \C^n$
which looks  similar to  the  Selberg integral and our integrals. For this reason,
it is natural to think that the story with beta-integrals of this type is not finished.
 \hfill $\boxtimes$

\sm

{\bf\punct Complex beta integrals and complex hypergeometric functions.}
Notice that the integral \eqref{eq:beta} for $\B^\C(\bfa,\bfb)$ is a twin of the classical Euler beta integral, this formula 
appeared in the book by Gelfand, Graev, Vilenkin \cite{GGV}, Section II.3.7, 1962, as an exercise (a careful proof is not completely
obvious). Our formulas \eqref{eq:main}, \eqref{eq:trapezoid-0}-\eqref{eq:trapezoid-1} also are twins of formulas from Neretin \cite{Ner-ray}.
It is well-known that the Selberg integral interpolates matrix beta integrals%
\footnote{On beta integrals, see Askey \cite{Ask}, Andrews, Askey, Roy \cite{AAR}, on matrix beta integrals, see a review  \cite{Ner-beta}. Beta integrals are integrals whose integrands are products
(for instance, of power functions or Gamma functions) and the right hand side is a product of Gamma-functions (or  $q$-Gamma function). Quite often such integrands are weight
functions for orthogonal polynomials or are related to explicitly solvable spectral problems with continuous or partially continuous spectra.}  of Hua Loo Keng type \cite{Hua} with respect to the dimension $2\theta$ of a basic field ($\R$, $\C$, or quaternions $\mathbb{H}$) and sizes of matrices.
One of purposes of \cite{Ner-ray} was an interpolation of  matrix beta integrals of more general types as Gindikin \cite{Gin}  and the author \cite{Ner-hua}. Another purpose was an interpolation of joint distributions of eigenvalues of growing
Hermitian matrices (for  stochastic processes generated 
by integrals from \cite{Ner-ray}, 
see Cuenca \cite{Cue}, Assiotis, Najnudel \cite{AN}). Our calculations in Section 2 are twins of calculations in \cite{Ner-ray}.

 Now there are known many
complex twins of real beta integrals,
see Ba\-zha\-nov, Mangazeev, Sergeev, \cite{BMS}, Kels \cite{Kel1}, \cite{Kel2}, Derkachov, Manashov \cite{DMV1}, \cite{DM2},  
Derkachov, Manashov, Valinevich \cite{DMV1}, \cite{DMV2}, Neretin \cite{Ner-dou}, Sarkissian, Spiridonov \cite{SS1}.
 Usually, these extensions are not
automatic, often usual ways are non-available  or  meet serious analytic difficulties. For instance, if we know the Selberg integral \eqref{eq:selberg}, then we can 
evaluate the integrals
\begin{equation*}
\int\limits_{\R^n} \prod_{k=1}^n x_k^{\sigma-1} e^{-ax_k}\! \prod_{1\le k<l\le n} \! |x_k-x_l|^{2\theta}\, dx
\,\, \,\text{and}\,\,  \int\limits_{\R^n} e^{-a\sum x_k^2} \prod_{1\le k<l\le n} |x_k-x_l|^{2\theta}\, dx
\end{equation*} 
by examination of asymptotics of \eqref{eq:selberg} as $\Re \tau\to\infty$ and $\Re\sigma$, $\Re\tau\to\infty$ respectively, see Andrews, Askey, Roy \cite{AAR}, Corollaries 8.2.2, 8.2.3. 
But we can not do the same procedure with their complex twins
\begin{equation*}
\int\limits_{\C^n} \prod_{k=1}^n z_j^{\sigmA-\1} e^{i\Re z_j} \prod_{1\le k<l\le n} |z_k-z_l|^{2\thetA}\, \dd{z}
 \,\,  \text{and}\,\,  \int\limits_{\C^n} e^{-i\Re\sum z_k^2} \prod_{1\le k<l\le n} |z_k-z_l|^{2\thetA}\, \dd{z}
\end{equation*} 
due to a narrow domain of convergence \eqref{eq:ds-convergence-1}--\eqref{eq:ds-convergence-3}. Also these integrals are not absolutely convergent, our way of evaluation of the
Dotsenko-Fateev  integral formally can be applied to these integrals, but this leads to too dangerous manipulations with non-absolutely convergent integrals  (cf. \cite{Ner-ray}, Subsect. 2.7, 2.8).

On the other hand, the Selberg integral has  different real forms, see \cite{AAR}, Chapter 8, Exercise 14, and
\cite{For}, formula (2.6).
Apparently, for the  complex case we have only one Dotsenko-Fateev  integral.

\sm

Another part of this story is hypergeometric functions of the complex field,  which are a kind of a lost chapter of theory of hypergeometric
functions%
\footnote{There are two groups of authors who worked on complex twins of hypergeometric functions and beta-integral,  mathematicians motivated by representations of the Lorentz
group and mathematical physicists, who also partially are motivated by representations of the Lorentz group.}, see Gelfand, Graev, Retakh \cite{GGR}, Ismagilov \cite{Ism},  Mimachi \cite{Mim}, Molchanov, Neretin \cite{MN}, Derkachov, Spiridonov,
\cite{DS}, Neretin \cite{Ner-bar}, Sarkissian, Spiridonov \cite{SS2}. For instance, the counterpart
of  the Gauss hypergeometric function $_2F_1$ is the following integral 
$$
{\vphantom {F^\C_1}}_2F_1^\C\left[\begin{matrix}\bfa,\bfb;\bfc \end{matrix};z \right]:=
\frac 1{\pi \B^\C(\bfa,\bfb)}\int_\C t^{\bfb-\1}(1-t)^{\bfc-\bfb-\1}(1-zt)^{-\bfa}\,\dd{t},
$$ 
defined in Gelfand, Graev, Retakh \cite{GGR} (this integral is a precise twin of the Euler integral representation
of the Gauss hypergeometric function $_2F_1$, see, e.g., \cite{AAR}, Sect.2.2). It admits an explicit evaluation as a quadratic expression
with two summands with products of Gaussian hypergeometric functions, see \cite{MN}, Theorem 3.9. Twins of generalized hypergeometric  functions $_pF_q$ were defined
in \cite{Ner-bar}. Apparently, the most of identities for usual hypergeometric functions $_pF_q$ have complex twins, as it is explained in \cite{Ner-bar}.

It is natural to ask about complex twins of the Heckman-Opdam multivariate hypergeometric functions \cite{HO}. This note gives a reason for this question,
 at least for 
functions related to the root systems of type $A_n$ (i.e., to the Jack functions).
Consider the space $\R^1\times \R^2\times \dots \times \R^n$ with coordinates $\cX=\{x_{j\alpha}\}$, where $1\le\alpha\le j\le n$. 
Consider  the  polyhedral cone $\cR_n\subset \prod_{j=1}^n \R^j$ consisting of $\cX$ satisfying the interlacing conditions
$$
x_{(j+1)1}\ge x_{j1} \ge x_{(j+1)2}\ge x_{j2}\ge x_{(j+1)3}\ge \dots\ge x_{(j+1)(j+1)}, \quad x_{nn}>0,
$$
see \cite{Ner-ray}, Subsect 1.1.
We fix 
$t_n\ge \dots\ge t_1>0$ and denote by $\cR_n(t)$ the set of $\cX\in \cR_n$ satisfying the conditions
$x_{n\alpha}=t_\alpha$
(i.e., we consider a section of the cone by an affine subspace). So we get a convex polyhedron in $\R^1\times\dots\times \R^{n-1}$.

According Okounkov and Olshanski \cite{OO} (see also Kazarnovski-Krol \cite{Kaz})
the Jack functions admit integral representations of the form
 \begin{multline}
 J_{\sigma_1,\dots,\sigma_n}^\theta(t_1,\dots,t_n)= C(\sigma,\theta)
 \prod_{\alpha=1}^n t_\alpha^{\sigma_n+(n-1)\theta/2}\prod_{1\le \alpha<\beta\le n} (t_\beta-t_\alpha)^{1-2\theta}
  \times\\ \times
 \int_{\cR_{\phantom{.}\,n}(t)} \prod_{j=1}^{n-1}\prod_{\alpha=1}^j x_{j\alpha}^{\sigma_j-\sigma_{j+1}-\theta}
   \times\\ \times
\prod_{j=1}^{n-1} \frac{\prod\limits_{1\le\alpha\le j;\,1\le p\le n}
|x_{j\alpha}-x_{(j+1)p}|^{\theta-1}\Bigr|_{x_{n1}=t_1, \dots, x_{nn}=t_n}}
{\prod\limits_{1\le\alpha<\beta\le j}|x_{j\alpha}-x_{j\beta}|^{2\theta-1}}
\prod_{1\le\alpha\le j\le n-1}dx_{j\alpha},
\label{eq:jack}
 \end{multline}
 where $C(\sigma,\theta)$ is a normalization constant (a product of Gamma-functions).
 For $\theta=1/2$, 1, 2 this formula gives spherical functions on the Riemannian symmetric spaces
\begin{equation}
 \GL(n,\R)/\mathrm{O}(n),\,\, \GL(n,\C)/\mathrm{U}(n),\,\, \GL(n,{\mathbb H})/\Sp(n),
 \label{eq:GL}
 \end{equation}
where $\mathbb{H}$ is the algebra of quaternions and $\Sp(n)$ is the quaternionic unitary group%
\footnote{These spaces are spaces of positive definite matrices of size $n$ over $\R$, $\C$,  $\mathbb{H}$.}. For the space  $\GL(n,\C)/{\mathrm U}(n)$,
  Gelfand and Naimark \cite{GN}, Sect. 9, obtained an elementary expression, this work was an initial point for further calculations of such type.

It is natural to hope that Jack functions of the complex field can be obtained from this expression by
the same substitution \eqref{eq:substitution}. It is natural to hope that they are eigenfunctions of  two commuting families of the Sekiguchi operators%
\footnote{See Sekiguchi \cite{Sek}, Macdonald \cite{Mac}, Sect. VII.3, Example 3.}
\begin{align*}
D(u;\thetA)=\prod_{1\le k<l\le n}(t_k-t_l)^{-1} \det\Bigl\{t_k^{n-l}\bigl(t_k\frac{\partial}{\partial t_k}+(n-l)\theta+u\bigr)\Bigr\}_{1\le k, l\le n};
\\
\ov D(u';\thetA)=\prod_{1\le k<l\le n}(\ov t_k-\ov t_l)^{-1} \det\Bigl\{\ov t_k^{n-l}\bigl(\ov t_k\frac{\partial}{\partial \ov t_k}+(n-l)\theta'+u'\bigr)\Bigr\}_{1\le k, l\le n}
\end{align*}
(cf.,  \cite{MN}, Propositions 3.9, 3.11, 3.12, \cite{Ner-bar},  Corollary 1.4). Other natural questions:
does exist a complex twin of the Harish-Chandra hypergeometric transform? Are spherical distributions on the symmetric spaces%
\footnote{These spaces are complexifications of the spaces \eqref{eq:GL}. They are spaces of nondegenerate complex matrices,   precisely,
  symmetric matrices of size $n$,
 square matrices of size $n$, skew-symmetric matrices of size $2n$. On spherical distributions on symmetric spaces, see Heckman, Schlichtkrull  \cite{HS}, Chapters 2, 3).} 
$$
 \GL(n,\C)/\mathrm{SO}(n,\C),\, \GL(n,\C)\times\GL(n,\C)/\mathrm{diag} \GL(n,\C),\, \GL(2n, \C)/\Sp(2n,\C)
$$
Jack functions of the complex field?.

\section{Evaluation of integrals}

\COUNTERS

{\bf \punct The Dirichlet integral.}
The complex twin of the  Dirichlet integral
(see, e.g. \cite{AAR}, Theorem I.8.6)
$$
\int\limits_{\substack{t_1>0,\dots,t_n>0\\ t_1+\dots+t_n<1}} \prod_{j=1}^{n} t_j^{a_j-1}\cdot \Bigl(1-\sum_{j=1}^n t_j\Bigr)^{a_{n+1}-1} \,dt_1\dots dt_n
=\frac{\prod_{j=1}^{n+1}\Gamma(a_j)}{\Gamma(\sum_{j=1}^{n+1} a_j\bigr)}
$$
 is given by the following statement.

\begin{proposition}
Denote
\begin{equation}
D_n^\C(\bfa_1,\dots, \bfa_{n+1}):=\int_{\C^n} \prod_{j=1}^{n} t_j^{\bfa_j-\1}\cdot \Bigl(1-\sum_{j=1}^n t_j\Bigr)^{\bfa_{n+1}-\1} \dd{t}_1\dots \dd{t}_n.
\label{eq:dirichlet0}
\end{equation}
Then
\begin{equation}
D_n^\C(\bfa_1,\dots, \bfa_{n+1})=
\pi^n \frac{\prod_{j=1}^{n+1}\Gamma^\C(\bfa_j)}{\Gamma^\C(\sum_{j=1}^{n+1} \bfa_j)}, 
\label{eq:dirichlet}
\end{equation}
the conditions of convergence are 
$$
\lfloor \bfa_j\rfloor>0, \qquad \lfloor\sum \bfa_j\rfloor<1.
$$
\end{proposition}

{\sc Proof.} Notice that for $u\in\C$ we have
\begin{equation}
\frac 1\pi \int_\C t^{\bfa-\1} (u-t)^{\bfb-\1}\,\dd{t}_1\dots \dd{t}_n= u^{\bfa+\bfb-\1}\,
\frac{\Gamma^\C(\bfa)\,\Gamma^\C(\bfb)} {\Gamma^\C(\bfa+\bfb)}.
\label{eq:beta-extended}
\end{equation}
Indeed, the substitution $t=uz$ reduces this integral to \eqref{eq:beta}.

Next, we integrate with respect to $t_n$ in \eqref{eq:dirichlet0}, for this purpose in 
\eqref{eq:beta-extended} we set $t=t_n$, $u=1-t_1-\dots-t_{n-1}$. We come to 
$$
D_n^\C(\bfa_1,\dots, \bfa_{n+1})=
\pi \, \frac{\Gamma^\C(\bfa_n)\,\Gamma^\C(\bfa_{n+1}) }{\Gamma^\C(\bfa_n+\bfa_{n+1})}\cdot
 D_{n-1}^\C(\bfa_1,\dots,\bfa_{n-1},\bfa_n+ \bfa_{n+1}).
$$
Iterating this identity, we get \eqref{eq:dirichlet}.
\hfill $\square$.

\sm

{\bf \punct The main lemma.}

\begin{lemma}
\label{l:main}
\begin{multline}
\int_{\C^n}
\prod_{p=1}^n
u_p^{\sigmA -\1} (1-u_p)^{\taU-\1} 
 \prod_{1\le\alpha\le n-1;\, 1\le p\le n} (u_p-z_\alpha)^{\thetA_\alpha-\1} \prod_{1\le p<q\le n}(u_p-u_q)^{\1}\,\dd{u}=
 \\=
(-1)^{\sum_{\alpha=1}^{n-1} \thetA_\alpha}\,\pi^n n! \frac{\Gamma^\C(\sigmA)\,\Gamma^\C(\taU)\prod_{\alpha=1}^{n-1}\Gamma^\C(\thetA_\alpha)}
{\Gamma^\C(\sigmA+\taU+\sum_{\alpha=1}^{n-1}\thetA_\alpha)}
\times\\\times
\prod_{\alpha=1}^{n-1} z_\alpha^{\sigmA+\thetA_\alpha-\1} (1-z_\alpha)^{\taU+\thetA_\alpha-\1}
\prod_{1\le \alpha,\beta\le n-1;\, \alpha\ne\beta}(z_\beta-z_\alpha)^{\thetA_\alpha-\frac{\1}{\2}}.
\label{eq:lemma-1}
\end{multline}
The domain of convergence is 
$$
\lfloor \sigmA\rfloor>0,\quad \lfloor \taU\rfloor>0, \quad \lfloor \thetA_\alpha\rfloor>0,
 \quad \lfloor\sigmA+\taU+\sum_{\alpha=1}^{n-1}\thetA_\alpha\rfloor<1.
$$
\end{lemma}

{\sc Remark.}
The last factor in the right hand side of \eqref{eq:lemma-1} can be represented in the form 
$$
\prod_{1\le \alpha,\beta\le n-1;\, \alpha\ne\beta} \!\!(z_\beta-z_\alpha)^{\thetA_\alpha} \prod_{1\le \alpha<\beta\le n-1}\!\!(z_\beta-z_\alpha)^{-\1}
=
\pm \prod_{1\le \alpha<\beta\le n-1}\!\!(z_\beta-z_\alpha)^{\thetA_\alpha+\thetA_{\beta}-\1}.
$$

{\sc Proof.} {\sc Step 1.} We define new variables $x_1$, \dots, $x_{n-1}$, $y\in \C$ instead of $u_p$
by 
\begin{align}
x_\alpha&=-\frac{\prod_{1\le p\le n}(u_p-z_\alpha)}
{\prod_{1\le \beta\le n-1, \beta\ne\alpha}(z_\beta-z_\alpha)};
\label{eq:change-1}
\\
y:&=\sum_{1\le p\le n} u_p- \sum_{1\le\beta\le n-1} z_\beta,
\label{eq:change-2}
\end{align}
this transformation is similar to Anderson \cite{And}.

Notice that any permutation of $u_p$ does not change $x_\alpha$, $y$. Let us show that this transformation is a bijection a.s.
between the quotient of $\C^n$ by the symmetric group $S_n$ and the space $\C^n$ of vectors $(x_1,\dots,x_{n-1},y)$.
For this purpose, we consider the following rational function:
$$
Q(t):=\frac{\prod_{1\le p\le n} (t-u_p)}{\prod_{1\le  \beta\le n-1}(t-z_\beta)}
$$
If $z_\alpha$ are pairwise different, then
$$
Q(t)
=t-\sum_{p=1}^{n} u_p+\sum_{\alpha=1}^{n-1} z_\alpha+
\sum_{\alpha=1}^{n-1} \frac{\res\limits_{t=z_\alpha} Q(t) }{t-z_\alpha}.
$$
The residues are
$
\res\limits_{t=z_\alpha}= -x_\alpha
$.
So
$$
Q(t)=t-y-\sum_{1\le \alpha\le n-1} \frac{x_\alpha}{t-z_\alpha},
$$
and $Q(t)$ is uniquely determined by $x_\alpha$ and $y$.

\sm

{\sc Step 2.}

\begin{lemma}
\label{eq:Jacobian1}
The complex Jacobian of the map \eqref{eq:change-1}-\eqref{eq:change-2}
is 
$$
J(u;z)=\frac{\prod_{1\le p\le q\le n}(u_p-u_q)}{\prod_{1\le\alpha<\beta\le n-1}(z_\beta-z_\alpha)}.
$$
\end{lemma}

The statement is Lemma 2.3.c from \cite{Ner-ray}.
It is a simple calculation,  which is reduced to the Cauchy determinant, see, e.g., \cite{Mac}, Sec. I.4, Example 6.

\sm

So the real Jacobian of the transformation is $|J(u;z)|^2=J(u;z)^\1$.

Denote by $I(u;z)$ the integrand in \eqref{eq:lemma-1}. Then our integral
is 
\begin{equation}
n!\int_{\C^n} I\bigl(u(x,y);z\bigr) J\bigl(u(x,y);z\bigr)^{-\1}\,\dd{x}\,\dd{y}.
\label{eq:I-J}
\end{equation}
The factor $n!$ arises since the map 
$(u_1,\dots,u_n)\mapsto (y,x_1,\dots,x_{n-1})$ is an $n!$-sheeted (ramified)
covering.
We must express our integrand in the variables $x$, $y$.

\sm

{\sc Step 3.}
Let us show that for any $a\in\C$ we have
\begin{equation}
\prod_{p=1}^n (u_p+a)=\Bigl( a+y-\sum_{\alpha=1}^{n-1} \frac{x_\alpha}{z_\alpha+a}\Bigr)
\prod_{\alpha=1}^{n-1} (z_\alpha+a).
\label{eq:sum-res}
\end{equation}
This is an identity `sum of residues is 0' for the rational function
$$R(t)=\frac{\prod_{1\le p\le n}(t-u_p)}{(t+a)\prod_{1\le\beta\le n-1}(t-z_\beta)}.
$$
Indeed,
the residues of $R(t)$ are
\begin{align*}
\res\limits_{t=z_\alpha}R(t)&=-\frac{x_\alpha}{z_\alpha+a},\qquad \res\limits_{t=\infty} R(t)=y+a,
\\
 \res\limits_{x=-a}R(t)&=
-\frac{\prod_{1\le p\le n}(u_p+a)}{\prod_{1\le\beta\le n-1}(z_\beta+a)}.
\end{align*}

\sm

{\sc Step 4.} Now we are ready to transform the integrand in \eqref{eq:I-J}. Setting
$a=0$, $-1$, $-z_\alpha$ in \eqref{eq:sum-res} we respectively  get 
\begin{align*}
\prod_p u_p^{\sigmA-\1}&=\prod_{\alpha=1}^{n-1} (z_\alpha)^{\sigmA-\1}
\Bigl( y-\sum_{\alpha=1}^{n-1} \frac{x_\alpha}{z_\alpha}\Bigr)^{\sigmA-\1};
\\
\prod_p (1- u_p)^{\taU-\1}&= (-1)^{n(\taU-\1)}\prod_{\alpha=1}^{n-1} (z_\alpha-1)^{\taU-\1}
\Bigl(-1+ y-\sum_{\alpha=1}^{n-1} \frac{x_\alpha}{z_\alpha-1}\Bigr)^{\taU-\1}=
\\
&=\prod_{\alpha=1}^{n-1} (1-z_\alpha)^{\taU-\1}
\Bigl(1- y-\sum_{\alpha=1}^{n-1} \frac{x_\alpha}{1-z_\alpha}\Bigr)^{\taU-\1};
\\
\prod_{p=1}^n (u_p-z_\alpha)^{\thetA_\alpha-\1}&=
\Bigl(-x_\alpha {\prod_{1\le \beta\le n-1, \beta\ne\alpha}(z_\beta-z_\alpha)}\Bigr)^{\thetA_\alpha-\1}.
\end{align*}
The factor $\prod(u_p-u_q)^\1$ appears in the numerator of $I(\cdot)$ and the denominator of $J(\cdot)^{-\1}$ and denominator, so it cancels.
So we come to the expression
\begin{equation*}
A\cdot \int_{\C^n} \Bigl(y-\sum_{1\le\alpha\le n-1}\frac{x_\alpha}{z_\alpha}\Bigr)^{\sigmA-\1}
\Bigl(1-y-\sum_{1\le\alpha\le n-1}\frac{x_\alpha}{1-z_\alpha}\Bigr)^{\taU-\1}
\prod_{\alpha=1}^{n-1} x_\alpha^{\thetA_\alpha-\1}\,\dd{x}\, \dd{y},
\end{equation*}
where
\begin{multline*}
A=(-1)^{\sum_{\alpha=1}^{n-1} \thetA_\alpha}\, n!
\times\\ \times
\prod_{\alpha=1}^{n-1} z_\alpha^{\sigmA-\1}(1-z_\alpha)^{\taU-\1} \!\! \prod_{1\le\alpha,\beta\le n-1,\,\alpha\ne\beta}
(z_\beta-z_\alpha)^{\thetA_\alpha-\1}\!\! \prod_{1\le\alpha<\beta\le n-1}
(z_\beta-z_\alpha)^{\1}
\end{multline*}
(the last two products are united in the final formula \eqref{eq:lemma-1}).

\sm

{\sc Step 5.} Changing the variables to 
$$
s=y-\sum_{\alpha=1}^{n-1} \frac{x_\alpha}{z_\alpha}  ,\qquad t_\alpha=\frac{x_\alpha}{z_\alpha(1-z_\alpha)}, 
$$
we come to the expression
\begin{equation*}
A\cdot \prod_{\alpha=1}^{n-1}(z_\alpha(1-z_\alpha))^{\thetA_\alpha} 
\int_{\C^n} s^{\sigmA-\1} \prod_{\alpha=1}^{n-1} t_\alpha^{\thetA_\alpha-\1}\cdot \Bigl(1-s-\sum_{\alpha=1}^{n-1} t_j\Bigr)^{\taU-\1}\,\dd{t}\,\dd{s}.
\end{equation*}
Applying the Dirichlet integral \eqref{eq:dirichlet} we get the desired statement.
\hfill $\square$

\sm

{\bf \punct Proof of Theorem \ref{th:1}.}
Denote
\begin{align*}
&\wh\sigmA_k:=\{\sigmA_j\}_{j=1,\dots, k}, 
\qquad
\wh \taU_k=\{\taU_j\}_{j=1,\dots, k},
\\
& \wh\thetA_k=\{\thetA_{j\alpha}\}_{j=1,\dots,k-1;\,\alpha=1,\dots,j}.
\end{align*}
Denote the integral \eqref{eq:main} by
$H(\wh\sigmA_n,\wh\taU_n,\wh\thetA_n)$. We integrate 
in the variables $z_{n1}$, \dots, $z_{nn}$ with Lemma \ref{l:main}
and come to
\begin{multline*}
H(\wh\sigmA_n,\wh\taU_n,\wh\thetA_n)=H(\wh\sigmA_{n-1},\wh\taU_{n-1},\wh\thetA_{n-1})
\times\\ \times
(-1)^{\sum_{\alpha=1}^{n-1}\thetA_{(n-1)\alpha}}\cdot \pi^n n! \cdot
 \frac{\Gamma^\C(\sigmA_n)\,\Gamma^\C(\taU_n)\prod_{\alpha=1}^{n-1}\Gamma^\C(\thetA_{(n-1)\alpha})}{\Gamma^\C(\sigmA_n+\taU_n+\sum_{\alpha=1}^{n-1}\thetA_{(n-1)\alpha})}.
\end{multline*}
Repeating the same operation we come to the desired statement.

\sm

{\bf \punct A dual lemma.} The following will be used to prove Lemma \ref{l:H}.

\begin{lemma}
\label{l:dual}
\begin{multline}
\int_{\C^n}\prod_{1\le\alpha\le n,\,1\le p\le n+1 }
(z_\alpha-u_p)^{\thetA-\1}\prod_{1\le \alpha<\beta\le n}
(z_\beta-z_\alpha)^\1\,\dd{z}=
\\=
n!\,\pi^n\,\frac{\Gamma^\C(\thetA)^{n+1}}{\Gamma^\C\bigl((n+1)\thetA\bigr)}\cdot \prod_{1\le p,q\le n+1;\, p\ne q} (u_p-u_q)^{\thetA-\frac\1 \2}
 .
\label{eq:int-dual}
\end{multline}
The condition of convergence of this integral is
$$
0<\lfloor \thetA\rfloor<\frac 1{n+1}.
$$
\end{lemma}

{\sc Proof.} Define new variables $w_1$, \dots, $w_{n+1}$ by 
$$
w_p:=\frac{\prod\limits_{1\le\alpha\le n}(z_\alpha-u_p)}{\prod\limits_{1\le q\le n+1,q\ne p}(u_q-u_p)},
$$
they are residues of the function
$$
H(t)=\frac{\prod\limits_{1\le\alpha\le n}(t-z_\alpha)}{\prod\limits_{1\le p\le n+1}(t-u_p)}
$$
at points $t=u_p$.
The residue at $\infty$ is $-1$, and therefore $\sum w_p=1$. On the other hand for fixed $u_p$
 a function $H(t)$ is uniquely determined by its residues, and therefore $w_p$ uniquely determine
 $z_\alpha$ up to a permutation.
 
  We choose $w_1$, \dots, $w_n$ as independent variables.
  
  \begin{lemma}
  \label{l:Jacobian2}
  The complex Jacobian of the transformation $(z,\dots,z_n)\mapsto(w_1,\dots, w_n)$ is
  $$
  J(z;u)=
  \det\limits_{1\le p,\,\alpha\le n}\Bigl\{\frac{\partial w_p}{\partial z_\alpha}\Bigr\}
= \frac{\prod\limits_{1\le\alpha<\beta\le n}(z_\beta-z_\alpha)}{\prod\limits_{1\le p<q\le n+1}(u_p-u_q)}.
  $$
  \end{lemma}
 
 {\sc Proof of Lemma \ref{l:Jacobian2}.}
 We have
 $$
 \frac{\partial w_p}{\partial z_\alpha}= \frac {w_p}{z_\alpha-u_p},
 $$
 So,
\begin{multline*}
J(z;u)=
%\frac{\cD(w_1,\dots,w_n)}{\cD(z_1,\dots,z_n)}=
\prod_{1\le p\le n} w_p\cdot
\det\limits_{1\le \alpha,\,p\le n}
\Bigl\{\frac 1{z_\alpha-u_p}\Bigr\}=
\\=
\frac{\prod\limits_{1\le\alpha,p\le n}(z_\alpha-u_p)}{\prod\limits_{1\le p, q\le n+1,q\ne p}(u_q-u_p)} \cdot
 \frac{\prod\limits_{1\le\alpha<\beta\le n}(z_\alpha-z_\beta)
  \prod\limits_{1\le p<q\le n}(u_p-u_q)}{\prod\limits_{1\le \alpha,p\le n}(z_\alpha-u_p)}.
\end{multline*}
Here we applied the formula for  the Cauchy determinant
$\det
\{\frac 1{z_\alpha-u_p}\}$ (see, e.g., \cite{Mac}, Sect I.4, Example 6).
After a cancellation we come to the desired formula.
\hfill $\square$

\sm 

The real Jacobian is $J(z;w)^{\1}$.

Next, we change variables $(z_1,\dots,z_n)\mapsto (w_1,\dots, w_n)$
in \eqref{eq:int-dual}. We have
\begin{align*}
(z_\beta-z_\alpha)^\1 \cdot  J(z;u)^{-\1}&=\prod_{1\le p<q\le n}(u_p-u_q)^\1;
\\
\prod_{1\le\alpha\le n,1\le p\le n+1 }
(z_\alpha-u_p)^{\thetA-\1}&=\prod_{p=1}^{n+1} w_p^{\thetA-\1}\cdot 
\prod_{1\le p<q\le n}(u_p-u_q)^{2\thetA-\2}.
\end{align*}

Now the integral in the left hand side of \eqref{eq:int-dual}   converts to
$$
n! \prod_{1\le p< q\le n+1}(u_p-u_q)^{2\thetA-\1} \int_{\C^n} \prod_{p=1}^{n+1} w_p^{\thetA-\1} \,\dd{w_1}\dots \dd{w_n}. 
$$
We remember that $w_{n+1}=1-\sum\limits_{p\le n} w_p$ and come to the Dirichlet integral \eqref{eq:dirichlet}.
This completes the proof of Lemma \ref{l:dual}.
%$\phantom{a}$ 
\hfill $\square$

\sm

{\bf\punct Proof of Theorem \ref{th:2}.}
Fix $m\le n$. Let $k< m$. Denote $\dd{\cZ}:=\prod_{j=k}^{n} \prod_{\alpha=1}^j \dd{z}_{j\alpha}$.
Denote
\begin{multline}
H_n^{m,k}:=
\int\limits_{\cC_{\,n}^{\,\vphantom{q_B}\,k}}%_{\cC_n^k}
\prod_{j=k}^{m-1}
\frac{\prod\limits_{1\le\alpha\le j;\, 1\le p\le j+1}(z_{j\alpha}-z_{(j+1)p})^{\nU-\1 }}
%{\prod\limits_{1\le \alpha<\beta\le j}(z_{j\alpha}-z_{j\beta})^{2\nU-\2}}
{\prod\limits_{1\le \alpha,\beta\le j,\,\alpha\ne\beta}(z_{j\alpha}-z_{j\beta})^{\nU-\1}}\times
 \\
 \times
 % \prod_{1\le \alpha<\beta\le k,\, \alpha\ne \beta}
 %(z_{\alpha k}-z_{\beta k})^{2\nU-\1}
 \prod_{1\le \alpha,\beta\le k,\, \alpha\ne \beta}
 (z_{k\alpha }-z_{k\beta })^{\nU-\frac\1\2}
\times\\ \times
 \prod_{j=m}^{n-1}\prod_{\alpha=1}^j z_{j\alpha}^{\sigmA_j-\sigmA_{j+1}-\thetA_{j\alpha}}(1-z_{j\alpha})^{\taU_j-\taU_{j+1}-\thetA_{j\alpha}}
 \prod_{p=1}^n z_{np}^{\sigmA_n-\1}(1-z_{np})^{\taU_n-\1}
 \times
 \\
 \times
 \prod_{j=m}^{n-1} \frac{\prod\limits_{1\le \alpha\le j;\, 1\le p\le j+1}(z_{j\alpha}-z_{(j+1)p})^{\thetA_{j\alpha}-\1}}
 {\prod\limits_{1\le\alpha,\beta\le j;\,\alpha\ne\beta}(z_{j\alpha}-z_{j\beta})^{\thetA_{j\alpha}-\1}}  \prod_{1\le p<q\le n}(z_{np}-z_{nq})^{\1}\,\,
 \dd{\cZ}.
 \label{eq:H}
 \end{multline}
 Since we have a cancellation in the first and second lines of the integrand,
 we can replace these rows  by
 \begin{multline}
 \prod\limits_{1\le\alpha\le k;\, 1\le p\le k+1}(z_{k\alpha}-z_{(k+1)p})^{\nU-\1}
 \cdot  \prod\limits_{1\le\alpha<\beta\le k}
 (z_{k\alpha}-z_{k\beta})^{\1 }
 \times
 \\
  \times
 \prod_{j=k+1}^{m-1}
\frac{\prod\limits_{1\le\alpha\le j;\, 1\le p\le j+1}(z_{j\alpha}-z_{(j+1)p})^{\nU-\1 }}
%{\prod\limits_{1\le \alpha<\beta\le j}(z_{j\alpha}-z_{j\beta})^{2\nU-\2}}
{\prod\limits_{1\le \alpha,\beta\le j,\,\alpha\ne\beta}(z_{j\alpha}-z_{j\beta})^{\nU-\1}}
.
\label{eq:two-rows}
\end{multline}

If $k=1$, then we have an integral over $\cC_n$ of the form
\eqref{eq:main}
with parameters 
\begin{align*}
&(\sigmA_m+(m-1)\nU,\sigmA_m+(m-2)\nU, \dots,\sigmA_m+2\nU,\sigmA_m+\nU, \sigmA_m, \sigmA_{m+1},\sigmA_{m+2}, \dots, \sigmA_n), 
\\
&(\taU_m+(m-1)\nU,\taU_m+(m-2)\nU, \dots,\taU_m+2\nU,\taU_m+\nU, \taU_m, \taU_{m+1},\taU_{m+2}, \dots, \taU_n), 
 \\
 &\bigl(\{\nU\}, \{\nU,\nU\},\dots,\underbrace{\{\nU,\dots\nU\}}_{\text{$m-1$ times}},
 \\
 & \qquad\qquad\qquad\{\thetA_{m\alpha}\}_{\alpha=1,\dots,m},\{\thetA_{(m+1)\alpha}\}_{\alpha=1,\dots,m+1},\dots, \{\thetA_{(n-1)\alpha}\}_{\alpha=1,\dots,n-1}\bigr).
 \end{align*}
Therefore, by Theorem \ref{th:1}
\begin{multline}
H_n^{m,1}=
(-1)^{m(m-1)\nU/2+\sum_{j=m}^{n-1}\sum_{\alpha=1}^j \thetA_{j\alpha}}\,
\pi^{n(n+1)/2}\, \prod_{j=1}^n j!
\times\\ \times 
 \prod_{j=1}^{m-1} \frac{\Gamma^\C(\sigmA_m+(m-j)\nU)\, \Gamma^\C(\taU_m+(m-j)\nU)}{\Gamma^\C(\sigmA_m+\taU_m+(2m-j-1)\nU)}
 \cdot
  \Gamma^\C(\nU)^{m(m-1)/2}  
    \times\\ \times
  \prod_{j=m}^{n-1}\prod_{\alpha=1}^j \Gamma^\C(\thetA_{j\alpha})
 \prod_{j=m}^n
 \frac{\Gamma^\C(\sigmA_j)\,\Gamma^\C(\taU_j)}{\Gamma^\C(\sigmA_j+\taU_j+\sum_{\alpha=1}^{j-1}\thetA_{(j-1)\alpha})}.
 \label{eq:Hnm1}
\end{multline}

\begin{lemma}
\label{l:H}
\begin{align*}
H_n^{m,k}&=H_n^{m,k+1}\cdot  \pi^k k! \frac{\Gamma^\C(\nU)^{k+1}}{\Gamma^\C((k+1)\nU)},\qquad \text{for $k<m-1$};
\\
H_n^{m,m-1}&=J_n^m(\wh \sigmA_n^m,\wh \taU_n^m;\nU, \wh \thetA_n^m)\cdot  \pi^{m-1} (m-1)!\,\, \frac{\Gamma^\C(\nU)^{m}}{\Gamma^\C(m\nU)}.
\end{align*}
\end{lemma}

{\sc Proof of Lemma \ref{l:H}.} Denote $\{z_{m\alpha}\}:=(z_{m1}, \dots, z_{mm})$. Fix $k$ to be $k\le m-1$.
The integrand in \eqref{eq:H} has the form
$$
A\bigl(\{z_{k\alpha}\}, \{z_{(k+1)\alpha}\}\bigr)\, 
B\bigl(\{z_{(k+1)\alpha}\}, \{z_{(k+2)\alpha}\}, \dots, \{z_{n\alpha\}}\bigr),$$
 where $A(\dots)$ is the first line of
\eqref{eq:two-rows} and $B(\dots)$ does not depend on
variables $z_{k1}, \dots, z_{kk}$ (it is the product of the second line
of \eqref{eq:two-rows} and the third and forth lines of
\eqref{eq:H}). So integration  with respect
to the variables $z_{k\alpha}$ gives
$$
\int_{\C^k} A\bigl(\{z_{k\alpha}\}, \{z_{(k+1)\alpha}\}\bigr)\,
\dd{z}_{k1} \dots \dd{z}_{kk}
\times B\bigl(\{z_{(k+1)\alpha}\}, \dots, \{z_{n\alpha}\}\bigr)
.
$$
We apply Lemma \ref{l:dual} and come to
$$
  \pi^k k! \frac{\Gamma^\C(\nU)^{k+1}}{\Gamma^\C((k+1)\nU)}
  \cdot \Bigl[ \prod_{\substack{1\le \alpha,\beta\le k+1\\ \alpha\ne \beta}}
 (z_{(k+1)\alpha }-z_{(k+1)\beta})^{\nU-\frac\1\2}
 \cdot B\bigl(\{z_{(k+1)\alpha}\}, \dots, \{z_{n\alpha}\}\bigr)
 \Bigr].
$$
The expression in the big square brackets is the integrand of $H_n^{m,k+1}$
in formula \eqref{eq:H}, when $k<m-1$ and the integrand in  $J_n^m(\wh\sigmA_n^m,\wh\taU_n^m;\nU,\wh\thetA_n^m)$
in formula \eqref{eq:trapezoid-0} when $k=m-1$. \hfill $\square$

\sm

Combining  Lemma \ref{l:H} with expression \eqref{eq:Hnm1} for $H_{n}^{m,1}$, we come
to  formula \eqref{eq:trapezoid-1}.  
The conditions of convergence of the $H_{n}^{m,1}$
are
\begin{align*}
&\lfloor\thetA_{j\alpha}\rfloor>0, \,\, \lfloor\taU_j\rfloor>0,\,\,  \lfloor \sigmA_j\rfloor>0, \,\,  \lfloor\sigmA_j+\taU_j+\sum_{\alpha=1}^{j-1}\thetA_{(j-1)\alpha} \rfloor<1,
\qquad \text{for $j> m$};
\\
&\lfloor\nU \rfloor>0, \,\,  \lfloor m\nU \rfloor<1,\,\, 
\lfloor \sigmA_m+\taU_m+2(m-1)\nU\rfloor<1,
\end{align*}
and these conditions are suffient for convergence of \eqref{eq:trapezoid-0}.

\sm

{\bf\punct The Dotsenko--Fateev--Aomoto integral.}
Let $m=n$. Then in \eqref{eq:trapezoid-0} we have integration 
$\dd{z}_{n1}, \dots, \dd{z}_{nn}$, two factors
$\prod_{j=m}^{n-1}$ disappear. The integrand transforms to
$$
\prod_{p=1}^n z_{np}^{\sigmA_n-\1}(1-z_{np})^{\taU_n-\1}
\cdot
\prod_{1\le\alpha,\beta\le n;\, \alpha\ne \beta} (z_{n\beta}-z_{n\alpha})^{\nU-\frac\1\2}
\cdot \prod_{1\le p<q\le n} (z_{np}-z_{nq})^\1.
$$
We unite the second and third factors and come to the desired form \eqref{eq:dfa}.

\tt

Fakult\"at f\"ur Mathematik, Universit\"at Wien;

Institute for Information Transmission Problems;

Moscow State University

e-mail: yurii.neretin(frog)univie.ac.at

URL: mat.univie.ac.at/$\sim$neretin

\end{document}